\documentclass[12pt]{article}
\usepackage{t1enc}
\usepackage[latin1]{inputenc}
\usepackage[english]{babel}
\usepackage{latexsym}

\usepackage{graphicx}
\usepackage{color}

\newcommand{\bp}{{\bf P}}

\newtheorem{conjecture}{Conjecture}[section]
\newtheorem{con}{Consequence}[section]

\newcommand{\bee}{{\bf e}}
\newcommand{\boo}{{\bf 0}}

\newcommand{\bs}{{\bf S}}

\newcommand{\bo}{{\bf 0}}

\newtheorem{theorem}{Theorem}[section]

\newtheorem{lemma}{Lemma}[section]

\def\kz{{\bf z}}

\def\bo{{\bf 0}}
\def\bee{{\bf e}}

\def\bs{{\bf S}}

\def\zd{{\cal Z}_d}
\def\z2{{\cal Z}_2}

\def\begg{\begin{equation}}
\def\endd{\end{equation}}
\def\bege{\begin{eqnarray}}
\def\ende{\end{eqnarray}}

\def\pe{{\bf P}}

\begin{document}

\centerline{\Large\bf Joint asymptotic behavior of local and
occupation times} \medskip \centerline{\Large\bf of random walk in higher
dimension }

\bigskip \bigskip \bigskip \bigskip \bigskip

\renewcommand{\thefootnote}{1} \noindent \textbf{Endre Cs\'{a}ki}\footnote{%
Research supported by the Hungarian National Foundation for Scientif\/ic
Research, Grant No. T 037886 and T 043037.}\newline
Alfr\'ed R\'enyi Institute of Mathematics, Hungarian Academy of Sciences,
Budapest, P.O.B. 127, H-1364, Hungary. E-mail address: csaki@renyi.hu

\bigskip

\renewcommand{\thefootnote}{2} \noindent \textbf{Ant\'{o}nia F\"{o}ldes}%
\footnote{%
Research supported by a PSC CUNY Grant, No. 66494-0035.}\newline
Department of Mathematics, College of Staten Island, CUNY, 2800 Victory
Blvd., Staten Island, New York 10314, U.S.A. E-mail address:
foldes@mail.csi.cuny.edu

\bigskip

\noindent \textbf{P\'al R\'ev\'esz}$^1$ \newline
Institut f\"ur Statistik und Wahrscheinlichkeitstheorie, Technische
Universit\"at Wien, Wiedner Hauptstrasse 8-10/107 A-1040 Vienna, Austria.
E-mail address: reveszp@renyi.hu

\bigskip \bigskip \bigskip

\noindent \textit{Abstract}: Considering a simple symmetric random
walk in dimension $d\geq 3$, we study the almost sure joint
asymptotic behavior of two objects:  first the local times of a
pair of neighboring points, then the local time of a point and the
occupation time of the surface of the unit ball around it.

\bigskip


\bigskip


\bigskip
\noindent Running title: Joint local and occupation time

\renewcommand{\thesection}{\arabic{section}.}

\section{Introduction and main results}

\renewcommand{\thesection}{\arabic{section}} \setcounter{equation}{0} %
\setcounter{theorem}{0} \setcounter{lemma}{0}

\noindent Consider a simple symmetric random walk $\{\mathbf{S}%
_n\}_{n=1}^{\infty} $ starting at the origin $\mathbf{0}$ on the $d$%
-dimensional integer lattice $\mathcal{Z}_d$, i.e. $\mathbf{S}_0=\mathbf{0}$%
, $\mathbf{S}_n=\sum_{k=1}^{n} \mathbf{X}_k$, $n=1,2,\dots$, where $\mathbf{X%
}_k,\, k=1,2,\dots$ are i.i.d. random variables with distribution
\[
\mathbf{P} (\mathbf{X}_1=\mathbf{e}_i)=\frac{1}{2d},\qquad
i=1,2,\ldots,2d
\]
and $\{\mathbf{e}_1,\mathbf{e}_2,...\mathbf{e}_d\}$ is a system of
orthogonal unit vectors in $\mathcal{Z}_d$ and $\mathbf{e}_{d+j}=-\mathbf{e}%
_j,$ $j=1,2,\ldots,d.$ Define the local time of the walk by
\begin{equation}
\xi(\mathbf{z},n):=\#\{k: \,\,0< k \leq n,\,\,\, \mathbf{S} _k=\mathbf{z}
\},\quad n=1,2,\ldots,  \label{loc1}
\end{equation}
where $\mathbf{z}$ is any lattice point of $\mathcal{Z}_d.$ The maximal
local time of the walk is defined as
\begin{equation}
\xi(n):=\max_{\mathbf{z} \in \mathcal{Z}_d}\xi(\mathbf{z},n).  \label{loc2}
\end{equation}

Define also the following quantities:
\begin{equation}
\eta(n):=\max_{0\leq k\leq n}\xi(\mathbf{S}_k,\infty),  \label{loc3}
\end{equation}
\begin{equation}
Q(k,n):=\#\{\mathbf{z}:\ \mathbf{z}\in\mathcal{Z}_d,\ \xi(\mathbf{z},n)=k\},
\end{equation}
\begin{eqnarray}
U(k,n)&:=& \#\{j:\ 0< j\leq n,\ \xi(\mathbf{S}_j,\infty)=k,\ \mathbf{S}%
_j\neq \mathbf{S}_\ell\ (\ell=1,2,\ldots,j-1)\}  \nonumber \\
&=&\#\{\mathbf{z}\in\mathcal{Z}_d:\, 0<\xi(\mathbf{z},n)\leq \xi(\mathbf{z}%
,\infty)=k\}.
\end{eqnarray}

Denote by $\gamma(n)=\gamma(n;d)$ the probability that in the first $n-1$
steps the $d$-dimensional path does not return to the origin. Then
\begin{equation}
1=\gamma(1)\geq \gamma(2)\geq ...\geq \gamma(n)\geq...>0.  \label{gain}
\end{equation}
It was proved in \cite{DE50} that

\textbf{Theorem A} (Dvoretzky and Erd\H os \cite{DE50}) \textit{For $d\geq 3
$
\begin{equation}
\lim_{n\to\infty}\gamma(n) =\gamma=\gamma(\infty;d)>0,  \label{gam1}
\end{equation}
and
\begin{equation}
{\gamma}<\gamma(n)<{\gamma}+O(n^{1-d/2}).  \label{gam2}
\end{equation}
Consequently
\begin{equation}
\mathbf{P}(\xi(\mathbf{0},n)=0,\, \xi(\mathbf{0},\infty)>0)=O\left(
n^{1-d/2}\right)
\end{equation}
as} $n\to\infty$.

So $\gamma$ is the probability that the $d$-dimensional simple symmetric
random walk never returns to its starting point.

Let $\xi(\mathbf{z},\infty)$ be the total local time at $\mathbf{z}$ of the
infinite path in $\mathcal{Z}_d$. Then for $d \geq 3$ (see Erd\H os and
Taylor \cite{ET60}) $\,\xi(\mathbf{0},\infty)$ has geometric distribution:

\begin{equation}
\mathbf{P}(\xi(\mathbf{0},\infty)=k)={\gamma}(1-{\gamma})^k,\qquad
k=0,1,2,...  \label{geo1}
\end{equation}
Erd\H{o}s and Taylor \cite{ET60} proved the following strong law for the
maximal local time:

\bigskip \noindent \textbf{Theorem B} (Erd\H os and Taylor \cite{ET60})
\textit{For $d\ge 3$
\begin{equation}
\lim_{n\to\infty}\frac{\xi(n)}{\log n}=\lambda \hspace{1cm} \mathrm{a.s.},
\label{la}
\end{equation}
where}
\begin{equation}
\lambda=\lambda_d=-\frac{1}{\log(1-{\gamma})}.
\end{equation}

\bigskip Following the proof of Erd\H{o}s and Taylor, without any new idea,
one can prove that
\begin{equation}
\lim_{n\to\infty}\frac{\eta(n)}{\log n}=\lambda \hspace{1cm} \mathrm{a.s.}
\label{lb}
\end{equation}

Erd\H{o}s and Taylor \cite{ET60} also investigated the properties of $Q(k,n)$%
. They proved

\bigskip\noindent \textbf{Theorem C} (Erd\H os and Taylor \cite{ET60})
\textit{For $d\geq 3$ and for any} $k=1,2,\ldots$
\begin{equation}
\lim_{n\rightarrow\infty}{\frac{Q(k,n)}{n}}= \gamma^2(1-\gamma)^{k-1}
\hspace{1cm} \mathrm{a.s.}  \label{qkn}
\end{equation}

Pitt \cite{Pitt} proved (\ref{qkn}) for general random walk and Hamana \cite
{Ham1}, \cite{Ham2} proved central limit theorems for $Q(k,n)$.

In \cite{CsFR05} we studied the question whether $k$ can be replaced by a
sequence $t(n)=t_n\nearrow\infty$ of positive integers in (\ref{qkn}). Let
\begin{equation}
\psi(n)=\psi(n,B)=\lambda\log n-\lambda B\log\log n.  \label{psi}
\end{equation}

\bigskip\noindent \textbf{Theorem D} \textit{\ Let $d\geq 3$, $%
\mu(t):=\gamma(1-\gamma)^{t-1}$ and $t_n:=[\psi(n,B)],\,\, (B>2)$, where
$\psi(n,B)$ is defined by {\rm (\ref{psi})}. Then we have
\[
\lim_{n\rightarrow\infty}\sup_{t\leq t_n}\left|{\frac{Q(t,n)}{n\gamma\mu(t)}}%
- 1\right|=\lim_{n\rightarrow\infty}\sup_{t\leq t_n}\left|{\frac{U(t,n)}{%
n\gamma\mu(t)}}- 1\right|=0\hspace{1cm} \mathrm{a.s.} \label{qtn}
\]
} Here in $\sup_{t\leq t_n}$, $t$ runs through positive integers.

For a set $A\subset\mathcal{Z}_d$ the occupation time of $A$ is defined by
\begin{equation}
\Xi(A,n):=\sum_{\mathbf{z}\in A}\xi(\mathbf{z},n).  \label{occup}
\end{equation}

Consider the translates of $A$, i.e. $A+\mathbf{u}=\{\mathbf{z}+\mathbf{u}%
:\, \mathbf{z}\in A\}$ with $\mathbf{u}\in\mathcal{Z}_d$ and define the
maximum occupation time by
\begin{equation}
\Xi^*(A,n):=\sup_{\mathbf{u}\in \mathcal{Z}_d}\Xi(A+\mathbf{u},n).
\label{maxop}
\end{equation}

It was shown in \cite{CsFRRS}

\bigskip\noindent \textbf{Theorem E} \textit{\ For $d\geq 3$ and for any
fixed finite set $A\subset\mathcal{Z}_d$
\begin{equation}
\lim_{n\to \infty}\frac{\Xi^*(A,n)}{\log n}=c_A\qquad\hspace{1cm} \mathrm{%
a.s.}  \label{lim}
\end{equation}
with some positive constant $c_A$, depending on $A$.}

\bigskip Now we present some more notations. For $\mathbf{z}\in \mathcal{Z}_d
$ let $T_\mathbf{z}$ be the first hitting time of $\mathbf{z}$, i.e. $T_%
\mathbf{z}:=\min\{i\geq 1:\mathbf{S}_i=\mathbf{z}\}$ with the convention
that $T_\mathbf{z}=\infty$ if there is no $i$ with
$\mathbf{S}_i=\mathbf{z}$. Let
$T=T_\bo$. In general, for a subset $A$ of $\mathcal{Z}_d$, let $T_A$ denote
the first time the random walk visits $A$, i.e. $T_A:=\min\{i\geq 1:\,
\mathbf{S}_i\in A\} =\min_{\mathbf{z}\in A}T_{\mathbf{z}}$. Let $\mathbf{P}_{%
\mathbf{z}}(\cdot)$ denote the probability of the event in the bracket under
the condition that the random walk starts from $\mathbf{z}\in\mathcal{Z}_d$.
We denote $\mathbf{P}(\cdot)= \mathbf{P}_{\mathbf{0}}(\cdot)$. Define
\begg
\gamma_{\mathbf{z}}:=\pe(T_{\mathbf{z}}=\infty).
\endd
Let $\mathcal{%
S}(r)$ be the surface of the ball of radius $r$ centered at the origin, i.e.
\[
\mathcal{S}(r):=\{\mathbf{z}\in\mathcal{Z}_d: \Vert \mathbf{z}\Vert=r\},
\]
where $\Vert\cdot\Vert$ is the Euclidean norm. Denote
\[
\Xi(\mathbf{z},n):=\Xi(\mathcal{S}(1)+\mathbf{z},n),
\]
i.e. the occupation time of the surface of the unit ball centered at
$\mathbf{z}\in\mathcal{%
Z}_d$.

Introduce further
\begin{equation}
p:=\mathbf{P}_{\mathbf{e}_1}(T_{\mathcal{S}(1)}<T).  \label{ap}
\end{equation}

In words, $p$ is the probability that the random walk, starting from $%
\mathbf{e}_1$ (or any other points of $\mathcal{S}(1)$), returns to $%
\mathcal{S}(1)$ before reaching $\mathbf{0}$ (including the case $T_{%
\mathcal{S}(1)}<T=\infty$).

In particular it was shown in \cite{CsFRRS}

\begin{equation}
\lim_{n\to \infty}\frac{\sup_{\kz\in\zd}\Xi(\kz,n)}{\log n}=
\frac1{-\log\left(p+\frac1{2d}\right)}=:\kappa\hspace{1cm} \mathrm{a.s.}
\label{lim1}
\end{equation}

\bigskip

It is easy to see that Theorem D implies

\begin{con}
With probability 1 there exists a random variable $n_0$ such that if $n\geq
n_0$ then for all $k=1,2\dots,\psi(n,B)$ there exist

\item  {\rm (i)} $\mathbf{z}\in\mathcal{Z}_d$ such that $\xi(\mathbf{z}%
,n)=k$,

\item  {\rm (ii)} $j\leq n$ such that $\xi(\mathbf{S}_j,\infty)=k$.
\end{con}

It would be interesting to investigate the joint behavior of the local time
of a point and the occupation time of a set, but in general this seems to
be a very complicated question so we will deal only with the following two
special cases. We will consider the joint behavior of the local times of
two neighboring points, and the local time of a point and the occupation
time of a ball of radius 1 centered at the point. Concerning the first
question one might like to know whether it is possible that in two
neighboring points the local times are simultaneously around $\lambda \log
n.$ More generally, we might ask whether the pairs of possible values of
\begin{equation}
(\xi(\mathbf{z},n),\, \xi(\mathbf{z}+\mathbf{e}_i,n))  \label{szom1}
\end{equation}
completely fill the lattice points in the set $(\log n)\mathcal{A}$  where
$\mathcal{A}$ is defined as
\[
\mathcal{A}:=\{(x,y)\in \mathcal{Z}_d:\, 0\leq x\leq \lambda,\, 0\leq
y\leq\lambda\}.
\]

The answer for this question turns out to be negative. However we will prove
that for
\begin{equation}
\mathcal{B}:=\{y\geq 0, x\geq 0:\, -(x+y)\log (y+x)+x\log x+y\log
y- (x+y)\log\alpha \leq 1\},  \label{defb}
\end{equation}
where
$$\alpha:=\frac{1-\gamma}{2-\gamma}$$
we have
\begin{theorem} Let $d\geq 4$.
For each $\varepsilon>0$, with probability 1 there exists an $%
n_0=n_0(\varepsilon)$ such that if $n\geq n_0$ then

\item  {\rm (i)} $(\xi(\mathbf{z},n),\xi(\mathbf{z}+\mathbf{e}_i,n))\in
((1+\varepsilon)\log n) \mathcal{B},\quad \quad\forall\mathbf{z}\in\mathcal{Z%
}_d, \quad\forall i=1,2,\ldots,2d $

\item  {\rm (ii)} for any $(k,\ell)\in ((1-\varepsilon)\log
n)\mathcal{B}%
\cap \mathcal{Z}_d$ and for arbitrary $i\in \{1,2,\ldots,2d\}$
there exist

random $\mathbf{z}_1,\mathbf{z}_2\in\mathcal{Z}_d$ for which
\[
(\xi(\mathbf{z}_1,n),\xi(\mathbf{z_1}+\mathbf{e}_i,n))=(k+1,\ell)
\]
\[
(\xi(\mathbf{z}_2,n),\xi(\mathbf{z_2}+\mathbf{e}_i,n))=(k,\ell+1).
\]
\end{theorem}

We will first show that without restriction on the dimension we have

\begin{theorem} Let $d\geq 3$.
For each $\varepsilon>0$, with probability 1 there exists an $%
n_0=n_0(\varepsilon)$ such that if $n\geq n_0$ then

\item  {\rm (i)}
$(\xi(\mathbf{S}_j,\infty),\xi(\mathbf{S}_j+\mathbf{e}_i,\infty))\in
((1+\varepsilon)\log n) \mathcal{B},\quad \quad\forall
j=0,1,2,\ldots,n,\quad\forall i=1,2,\ldots,2d$

\item  {\rm (ii)} for any $(k,\ell)\in ((1-\varepsilon)\log
n)\mathcal{B}%
\cap \mathcal{Z}_d$ and for arbitrary $i\in\{1,2,\ldots,2d\}$
there exists a

random integer $j=j(k,\ell)\leq n$ for which
\[
(\xi(\mathbf{S}_j,\infty),\xi(\mathbf{S}_j+\mathbf{e}_i,\infty))=(k+1,\ell).
\]
\end{theorem}

Concerning the occupation time of the unit ball, Consequence 1.1 and
Theorem E suggest the following

\begin{conjecture}
For any $\varepsilon>0$ with probability 1 there exists a random
variable $ n_0=n_0(\varepsilon)$ such that if $n\geq n_0$ then
for all $k=1,2,\dots, [(1-\varepsilon)\kappa\log n]$ there exists
$\mathbf{z}\in\mathcal{Z}_d$ such that $\Xi(\mathbf{z},n)=k.$
\end{conjecture}

A simple consequence of our Theorem 1.3 is that Conjecture 1.1 is true. As
we indicated above, we are interested in the joint asymptotic behavior of
the random sequence
\[
(\xi(\mathbf{z},n),\Xi(\mathbf{z},n)),\quad \mathbf{z}\in\mathcal{Z}_d
\]
as $n\to\infty$. One might ask again whether this random vector will fill
out all the lattice points of the triangle $(\log n)\mathcal{C}$, where
\[
\mathcal{C}=\{(x,y)\in \mathcal{Z}_d:\, 0\leq x\leq \lambda,\, x\leq y\leq
\kappa\}.
\]
As before, it turns out that the above triangle will not be filled.
Instead, we will prove the following theorem.

Define the set $\mathcal{D}$ as
\begin{equation}
\mathcal{D}:=\{y\geq x\geq 0:\, -y\log y+x\log
(2dx)+(y-x)\log((y-x)/p)\leq 1\}, \label{defd}
\end{equation}
where $p$ was defined in (\ref{ap}) and its value in terms of $\gamma$
is given by (\ref{pe}) below.

\begin{theorem} Let $d\geq 4$.
For each $\varepsilon>0$ with probability 1 there exists an $%
n_0=n_0(\varepsilon)$ such that if $n\geq n_0$ then

\item {\rm (i)} $(\xi(\mathbf{z},n),\Xi(\mathbf{z},n))\in
((1+\varepsilon)\log n) \mathcal{D},\quad \quad\forall\mathbf{z}\in\mathcal{Z%
}_d$

\item {\rm (ii)} for any $(k,\ell)\in ((1-\varepsilon)\log
n)\mathcal{D}%
\cap \mathcal{Z}_d$ there exists a random $\mathbf{z}\in\mathcal{Z}_d$ for
which
\[
(\xi(\mathbf{z},n),\Xi(\mathbf{z},n))=(k,\ell+1).
\]
\end{theorem}

\begin{theorem} Let $d\geq 3$.
For each $\varepsilon>0$ with probability 1 there exists an $%
n_0=n_0(\varepsilon)$ such that if $n\geq n_0$ then

\item  {\rm (i)}
$(\xi(\mathbf{S}_j+\bee_i,\infty),\Xi(\mathbf{S}_j+\bee_i,\infty)\in
((1+\varepsilon)\log n) \mathcal{D},\quad \forall j=1,2,\ldots,n,\quad
\forall i=1,2,\ldots,2d$

\item  {\rm (ii)} for any $(k,\ell)\in ((1-\varepsilon)\log
n)\mathcal{D}%
\cap \mathcal{Z}_d$ and for arbitrary $i\in\{1,2,\ldots,2d\}$ there
exists a

random integer $j=j(k,\ell)\leq n$ for which
\[
(\xi(\mathbf{S}_j+\bee_i,\infty),\Xi(\mathbf{S}_j+\bee_i,\infty))=(k,\ell+1).
\]
\end{theorem}

\bigskip\noindent
{\bf Remark 1.1} The condition $d\geq 4$ in Theorem 1.1 and
Theorem 1.3 is needed only for the convergence of (\ref{conv})
while proving parts (i). The proofs of parts (ii) in both theorems
work also for $d=3$.

\renewcommand{\thesection}{\arabic{section}.}

\section{Preliminary facts and results}

\renewcommand{\thesection}{\arabic{section}} \setcounter{equation}{0}
\setcounter{theorem}{0} \setcounter{lemma}{0}

Recall the definition of $\gamma$, $\gamma_\kz$, $T$ and $T_{\mathbf{z}}$
in Section 1.
\begin{lemma}For $ i=1, 2,\dots, 2d$
\begin{eqnarray}
\gamma_{\bee_i}&=&\gamma,
\label{same1}\\
\mathbf{P}(T<T_{\mathbf{e}_i})&=&\mathbf{P}(T_{\mathbf{e}_i}<T)=
\frac{1-\gamma}{2-\gamma}=\alpha,
\label{same2}\\
\mathbf{P}(T=T_{\mathbf{e}_i}=\infty)&=&\frac{\gamma}{2-\gamma}
=1-2\alpha.
\label{same3}
\end{eqnarray}
\end{lemma}

\vspace{2ex}\noindent \textbf{Proof.} By symmetry
$\gamma_{\bee_i}=\gamma_{\bee_1 },\,\, i=1, 2,\dots, 2d.$  Hence

$$1-\gamma=\sum_{i=1}^{2d}\pe(\bs_1=\bee_i)(1-\gamma_{\bee_i})=
\sum_{i=1}^{2d}\frac{1}{2d}(1-\gamma_{\bee_i})=(1-\gamma_{\bee_1})$$
thus we have (\ref{same1}). Furthermore  observe that

$$1-\gamma=\mathbf{P}(T<\infty)=
\mathbf{P}(T<T_{\mathbf{e_i}})+
\mathbf{P}(T_{\mathbf{e_i}}<T)\pe_{\bee_i}(T<\infty)$$
and
$$1-\gamma=\mathbf{P}(T_{\bee_i}<\infty)=\mathbf{P}(T_{\bee_i}<T)+
\mathbf{P}(T<T_{\bee_i})\mathbf{P}(T_{\bee_i}<\infty).$$
Solving this system of equations for $\mathbf{P}(T_{\mathbf{e_i}}<T)$
and $\mathbf{P}(T<T_{\bee_i})$, we get (\ref{same2}), and (\ref{same3})
follows from $\mathbf{P}(T=T_{\mathbf{e}_i}=\infty)=1-
\mathbf{P}(T<T_{\mathbf{e}_i})-\mathbf{P}(T_{\mathbf{e}_i}<T)$.

\begin{lemma}
For $i=1,2,\ldots, 2d$
\begin{equation}
\mathbf{P}(\xi(\mathbf{0},\infty)=k,
\xi(\mathbf{e}_i,\infty)=\ell)=(1-2\alpha) {{k+\ell
}\choose{k}}\alpha^{k+\ell}, \quad
k,\ell=0,1,\ldots  \label{gen1}
\end{equation}
\end{lemma}
\vspace{2ex}\noindent \textbf{Proof.} By (\ref{same2}), the
probability of $k$ visits in $\mathbf{0}$ and $\ell$ visits in
$\mathbf{e}_i$ in any particular order is $\alpha^{k+\ell}$. The
binomial coefficient in (\ref{gen1}) is the number of possible
orders. Finally, observe that  starting from either of the two
points, the probability that the walk does not return back to the
starting point, nor to the other point is $1-2\alpha.$
 Hence the lemma follows. $\Box$

Recall the definition of $p$ in (\ref{ap}).

\begin{lemma}
\begin{eqnarray}
p&=&1-\frac1{2d(1-\gamma)},  \label{pe} \\
\mathbf{P}(\Xi(\boo,\infty)=j)&=&
\left(1-p-\frac1{2d}\right)\left(p+\frac1{2d}\right)^{j-1},\quad j=1,2\ldots,
\label{ageo2} \\
\mathbf{P}(\xi(\boo,\infty)=k, \, \Xi(\mathbf{0},\infty)=\ell+1 )&=&
{{
\ell }\choose{k}}
\left(1-p-\frac1{2d}\right)p^{\ell-k}\left(\frac1{2d}\right)^k,
\label{bgeo3} \\
&&\ell=0,1,\ldots,\,k=0,1,\ldots,\ell.  \nonumber
\end{eqnarray}
\end{lemma}

\vspace{2ex}\noindent \textbf{Proof.} Let $Z(A)$ denote the number of
visits in the set $A$ up to the first return to zero. Clearly
\begin{equation}
\mathbf{P}(Z(\mathcal{S}(1))=j,T<\infty)=p^{j-1}\frac{1}{2d},\quad
j=1,2,\dots  \label{zs1}
\end{equation}
Summing up (\ref{zs1}) in $j$, we get
\begin{equation}
1-\gamma=\mathbf{P}(T<\infty)=\sum_{j=1}^{\infty}p^{j-1}\frac{1}{2d} =\frac{1%
}{2d(1-p)},
\end{equation}
implying (\ref{pe}).

Introduce further
\[
\tau =\sum_{j=1}^{\infty}I\{\mathbf{S}_{j}\in \mathcal{S}(1),\Vert\mathbf{S}%
_{j+1}\Vert>1\},
\]
thus $\tau$ is the number of outward excursions from $\mathcal{S}(1)$ to $%
\mathcal{S}(1)$, including the last incomplete one. Hence
\[
\Xi(\mathbf{0}, \infty)=\tau+\xi(\mathbf{0},\infty).
\]
Since $p$ is the probability that the random walk starting from any point of
$\mathcal{S}(1)$ returns to $\mathcal{S}(1)$ from outside, while $1/(2d)$ is
the probability of the same return through the origin, $p+1/(2d)$ is the
probability that the random walk, starting from any point of $\mathcal{S}(1)$%
, returns to $\mathcal{S}(1)$ in finite time, (\ref{ageo2}) is immediate.
Furthermore, it is easy to see that

\[
\mathbf{P}(\xi(\boo,\infty)=k, \, \tau=M+1 )= {{k+M }\choose{k}}
\left(1-p-\frac1{2d}\right)p^{M}\left(\frac1{2d}\right)^k,
\]
implying (\ref{bgeo3}). $\Box$

Recall and define
\begin{eqnarray}
\gamma_{\mathbf{z}}&:=&\mathbf{P}(T_{\mathbf{z}}=\infty),
\qquad \gamma_{\mathbf{z}}(n):=\mathbf{P}(T_{\mathbf{z}}\geq n),\\
q_{\mathbf{z}}&:=&\mathbf{P}(T<T_{\mathbf{z}}),\qquad
q_\kz(n):=\mathbf{P}(T<\min(n,T_\kz)),
\label{aq} \\
s_{\mathbf{x}}&:=&\mathbf{P}(T_{\mathbf{x}}<T),\qquad
s_\kz(n):=\mathbf{P}(T_\kz<\min(n,T)).
\label{as}
\end{eqnarray}

Moreover, put
$$
p(n):=\mathbf{P}_{\mathbf{e}_1}(T_{S(1)}<\min(n,T)).
$$

Similarly to Theorem A, we prove

\begin{lemma}
\begin{eqnarray}
1-{\gamma}_\mathbf{z}+\frac{O(1)}{n^{d/2-1}}&\leq&
1-\gamma_\mathbf{z}(n)\leq 1-\gamma_\mathbf{z},  \label{gn}\\
q_\kz+\frac{O(1)}{n^{d/2-1}}&\leq&
q_\kz(n)\leq q_\kz, \label{qn} \\
s_\kz+\frac{O(1)}{n^{d/2-1}}&\leq&
s_\kz(n)\leq s_\kz \label{sn} \\
p+\frac{O(1)}{n^{d/2-1}}&\leq&
p(n)\leq p,\label{pn}
\end{eqnarray}
and $O(1)$ is uniform in $\mathbf{z}$.
\end{lemma}

\noindent \textbf{Proof.} For the proof of (\ref{gn}) see Jain and Pruitt
\cite{JP}.

To prove (\ref{qn}) and (\ref{sn}), observe that
\begin{eqnarray*}
0\leq q_\kz-q_\kz(n)&=&\mathbf{P}(T<T_\kz,\,n\leq T<\infty) \leq
\mathbf{P}(n\leq T<\infty)=\gamma(n)-\gamma, \\
0\leq s_\kz-s_\kz(n)&=&\mathbf{P}(T_\kz< T,\,n\leq T_\kz <\infty) \leq
\mathbf{P}(n\leq T_\kz<\infty)=\gamma_\kz(n)-\gamma_\kz.
\end{eqnarray*}

To prove (\ref{pn}), introduce
$\mathbf{b}_j=\mathbf{e}_1+\mathbf{e}_j,\,\,$
$j=1,2,\ldots,2d,$ then we have
\begg
0\leq p-p(n)=\mathbf{P}_{\mathbf{e}_1}(n\leq
T_{\mathcal{S}(1)}<\infty)
=\sum_{j=1}^{2d} \mathbf{P}_{\mathbf{e}_1}(\mathbf{S}_1=\mathbf{b}_j,
n\leq T_{\mathcal{S}(1)}<\infty).  \label{p--}
\endd
Observe that by (\ref{gn}), each term in the above sum can be estimated
by
\[
\mathbf{P}_{\mathbf{e}_1}(\mathbf{S}_1=\mathbf{b}_j, n\leq T_{\mathcal{%
S}(1)}<\infty)=\frac{1}{2d}\mathbf{P}_{\mathbf{b}_j}(n-1\leq T_{\mathcal{%
S}(1)}<\infty)=\frac{O(1)}{n^{d/2-1}},
\]
proving the lemma. $\Box$

\begin{lemma}
For $i=1,2,\ldots,2d$, $k+\ell>0$, $n>0$ we have
\begg
\bp(\xi(\boo,n)=k,\xi(\bee_i,n)=\ell)\leq {k+\ell\choose k}
\alpha^{k+\ell},
\label{aln}
\endd
and for $i=1,2,\ldots,2d$, $\ell>0$, $n>0$ we have
\begg
\mathbf{P}(\xi(\bee_i,n)=k, \, \Xi(\mathbf{e}_i,n)=\ell)\leq
{\ell\choose{k}} p^{\ell-k}\left(\frac1{2d}\right)^k.
\label{pnn}
\endd
\end{lemma}

\noindent \textbf{Proof.} To show (\ref{aln}), recall that by Lemma 2.1,
$q_{\bee_i}=s_{\bee_i}=\alpha$. The time between consecutive visits to
$\boo$ or $\bee_i$ is less than $n$, hence using the upper inequalities
in (\ref{qn}) and (\ref{sn}), it is easy to see that the probability of
$k$ visits in $\boo$ and $\ell$ visits in $\bee_i$ up to time $n$ in any
particular order, is less than $\alpha^{k+\ell}$. Now (\ref{aln}) is
seen by observing that the number of particular orders is the binomial
coefficient in (\ref{aln}).

Similarly, we can get (\ref{pnn}) by using (\ref{pn}). $\Box$

\renewcommand{\thesection}{\arabic{section}.}
\section{The basic equations}

\renewcommand{\thesection}{\arabic{section}} \setcounter{equation}{0} %
\setcounter{theorem}{0} \setcounter{lemma}{0}

It follows from Lemma 2.2 and Stirling formula that the asymptotic
relation
\begin{equation}
\log\mathbf{P}(\xi(\boo,\infty)=[x\log n],\,
\xi(\bee_i,\infty)=[y\log n]) \sim -g(x,y)\log n,\quad n\to\infty
\label{joint1}
\end{equation}
 holds for $i\in\{1,2,\ldots,2d\}$, $x \geq 0, \,y\geq 0$, where
\[
g(x,y)=-(x+y)\log (y+x)+x\log x+y\log y- (x+y)\log\alpha.
\]

It follows that $\mathbf{P}(\xi(\boo,\infty)=[x\log n],\,
\xi(\bee_i,\infty)=[y\log n])$ is of order $1/n$ if $(x,y)$ satisfies
the basic equation

\begg g(x,y)=1,\quad x\geq 0, \,y\geq 0. \label{basic1}
\endd
Observe that $g(x,y)$ is the function defining the set ${\cal
B}$  in (\ref{defb}). The next lemma describes the major properties of
the boundary of the set ${\cal B}$.

\begin{lemma}
\item {\rm (i)} For the points $(x,y)$ satisfying {\rm (\ref{basic1})}
we have
\begin{eqnarray}
x_{max}=y_{max}&=&\lambda,\\
(x+y)_{max}&=&\frac{1}{\log{\frac{1}{2\alpha}}},
\end{eqnarray}

when this maximum occurs then $x=y.$

\item {\rm (ii)} If $x=x_{max}=\lambda$, then $y=\lambda(1-\gamma)$ and
vica versa.

If $x=0$, then $y=\frac{1}{\log(1/\alpha)}$ and vica versa.

\item {\rm (iii)} For a given $x$, the equation {\rm
(\ref{basic1})} has one solution in $y$ for $x<x_0$, and for
$x=\lambda$ and

two solutions for $x_0\leq x<\lambda$, where
$$x_0=\frac{1}{\log (1/\alpha)}.$$
\end{lemma}

\begin{figure}[htbp]
  \centering
\begin{picture}(0,0)%
\includegraphics{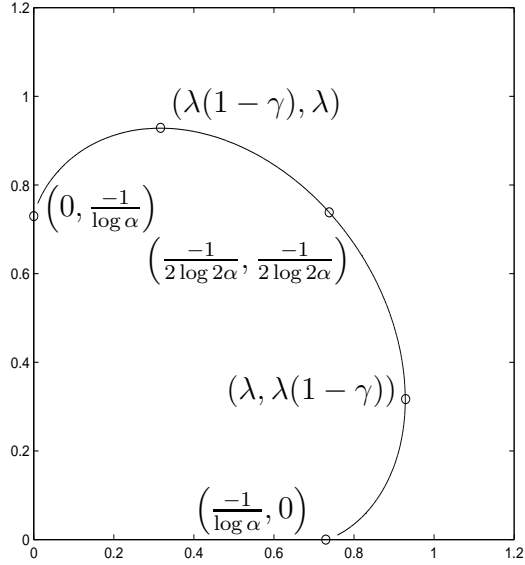}%
\end{picture}%
\setlength{\unitlength}{3947sp}%
\begingroup\makeatletter\ifx\SetFigFont\undefined%
\gdef\SetFigFont#1#2#3#4#5{%
  \reset@font\fontsize{#1}{#2pt}%
  \fontfamily{#3}\fontseries{#4}\fontshape{#5}%
  \selectfont}%
\fi\endgroup%
\begin{picture}(3319,3600)(301,-2761)
\put(1351, 91){\makebox(0,0)[lb]{\smash{{\SetFigFont{12}{14.4}{\familydefault}{\mddefault}{\updefault}{\color[rgb]{0,0,0}$(\lambda(1-\gamma),\lambda)$}%
}}}}
\put(1501,-2461){\makebox(0,0)[lb]{\smash{{\SetFigFont{12}{14.4}{\familydefault}{\mddefault}{\updefault}{\color[rgb]{0,0,0}$\left(\frac{-1}{\log\alpha},0\right)$}%
}}}}
\put(1191,-886){\makebox(0,0)[lb]{\smash{{\SetFigFont{12}{14.4}{\familydefault}{\mddefault}{\updefault}{\color[rgb]{0,0,0}$\left(\frac{-1}{2\log2\alpha},\frac{-1}{2\log2\alpha}\right)$}%
}}}}
\put(1700,-1711){\makebox(0,0)[lb]{\smash{{\SetFigFont{12}{14.4}{\familydefault}{\mddefault}{\updefault}{\color[rgb]{0,0,0}$(\lambda,\lambda(1-\gamma))$}%
}}}}
\put(551,-561){\makebox(0,0)[lb]{\smash{{\SetFigFont{12}{14.4}{\familydefault}{\mddefault}{\updefault}{\color[rgb]{0,0,0}$\left(0,\frac{-1}{\log\alpha}\right)$}%
}}}}
\end{picture}%

  \caption{The set $\mathcal{B}$ in the case of the two-point set, $d=3$.}
  \label{fig:twop}
\end{figure}

\bigskip\noindent \textbf{Proof.} Differentiating (\ref{basic1}) as
an implicit function of $x,$  $y$ takes its maximum ($y'=0$) at
$x=\lambda(1-\gamma)$ and the value of this maximum is
$y=\lambda,$  which proves the first statements in (i) and (ii).

Similarly, if we maximize $x+y$ as a function of $x$ (i.e.
$1+y'=0$) then we get that this occurs when $x=y$ and the second
part of (i) follows.

Solving (\ref{basic1}) when $x=0$ for $y$, we get the second part
of (ii).

Now we turn to the proof of (iii). For given $0\leq x\leq \lambda$
consider $g(x,y)$ as a function of $y$. We have
\[
{\frac{\partial g}{\partial y}}=\log{\frac{y}{\alpha(x+y)}}
\]
and this is equal to zero if $y=x(1-\gamma)$. It is easy to see
that $g$ takes a minimum here and is decreasing in $(0,x(1-\gamma))$ and
increasing in $(x(1-\gamma),\lambda)$. Moreover,
\[
{\frac{\partial^2 g}{\partial y^2}}={\frac{1}y}-{\frac{1}{x+y}}>0,
\]
hence $g$ is convex from below. We have for $0\leq x<\lambda,$
that this minimum is
\[
g(x,x(1-\gamma))=\frac{x}{\lambda}<1,
\]
and one can easily see that
\[
g(x,0)=x\log(1/\alpha)\left\{
\begin{array}{ll}
& < 1\,\, \mathrm{if}\,\, x<x_0, \\
& = 1\,\, \mathrm{if}\,\, x=x_0, \\
& > 1\,\, \mathrm{if}\,\, x>x_0.
\end{array}
\right.
\]
This shows that equation (\ref{basic1}) has one solution if
$0\leq x<x_0$ and two solutions if $x_0\leq x<\lambda$.

For $x=\lambda$, it can be seen that $y=\lambda (1-\gamma)$ is the
only solution of $g(x,y)=1$.

The proof of Lemma 3.1 is complete. $\Box$

For further reference introduce the following notations to
describe the boundary of ${\cal B}$: for $x_0\leq x<\lambda$ let
$y_{1,{\cal B}}(x)<y_{2,{\cal B}}(x)$ denote the two solutions and
for $0\leq x<x_0$ let $y_{2,{\cal B}}(x)$ denote the only solution
of (\ref{basic1}). Define $y_{1,{\cal B}}(x)=0$ for $0\leq x<x_0$
and $y_{1,{\cal B}}(\lambda)=y_{2,{\cal
B}}(\lambda)=\lambda(1-\gamma)$. Then the set ${\cal B}$ can be
given as
$$
{\cal B}=\{0\leq x\leq \lambda,\, y_{1,{\cal B}}(x)\leq y\leq
y_{2,{\cal B}}(x)\}.
$$

For further discussion of properties of the set ${\cal B}$ see Section
6.

Concerning similar description of the set ${\cal D}$ belonging to
the other problem, it follows from (\ref{bgeo3}) of Lemma 2.3 and
Stirling formula that the asymptotic relation
\begin{equation}
\log\mathbf{P}(\xi(\boo,\infty)=[x\log n],\, \Xi(\boo,\infty)=[y\log n])
\sim
-f(x,y)\log n,\quad n\to\infty  \label{joint2}
\end{equation}
holds for $0\leq x\leq y$, where
\[
f(x,y)=-y\log y+x\log x+(y-x)\log(y-x)+x\log(2d)+(y-x)\log(1/p).
\]

It follows that $\mathbf{P}(\xi(\boo,\infty)=[x\log n],\,
\Xi(\boo,\infty)=[y\log n])$ is of order $1/n$ if $(x,y)$ satisfies the
basic equation
\begin{equation}
f(x,y)=1,\qquad 0\leq x\leq y  \label{basic2}.
\end{equation}

\begin{lemma}
\item {\rm (i)} For the maximum values of $x,y$, satisfying {\rm
(\ref{basic2})}, we have
\begin{eqnarray}
x_{\max}&=&{\frac{1}{\log(2d(1-p))}}=\lambda,\label{xmax} \\
y_{\max}&=&{\frac{-1}{\log(p+\frac1{2d})}}=\kappa.
\end{eqnarray}
\item{\rm (ii)} If $x=x_{\max}=\lambda$, then $y=\lambda/(1-p)$. If
$y=y_{\max}=\kappa$, then $x=\kappa/(2dp+1)$. If $x=0$,

then $y=1/\log(1/p)$.

\item {\rm (iii)} For given $x$ the equation {\rm (\ref{basic2})}
has one solution in $y$ for $0\leq x< 1/\log(2d)$ and for

$x=\lambda$, and two solutions in $y$ for $\, 1/\log(2d)\leq
x<\lambda$.
\end{lemma}
\begin{figure}[htbp]
  \centering
\begin{picture}(0,0)%
\includegraphics{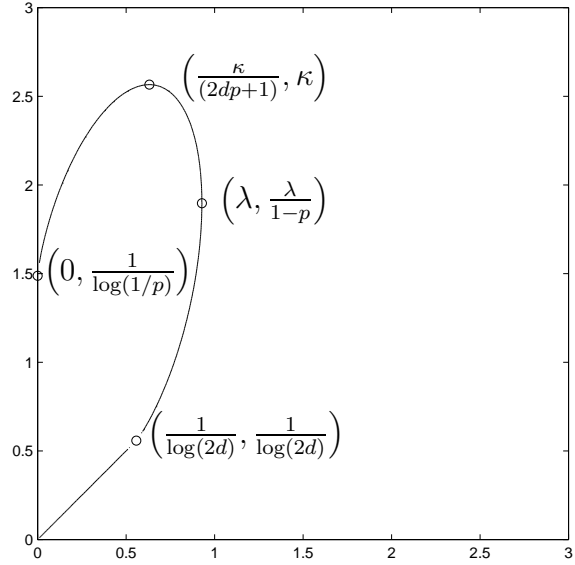}%
\end{picture}%
\setlength{\unitlength}{3947sp}%
\begingroup\makeatletter\ifx\SetFigFont\undefined%
\gdef\SetFigFont#1#2#3#4#5{%
  \reset@font\fontsize{#1}{#2pt}%
  \fontfamily{#3}\fontseries{#4}\fontshape{#5}%
  \selectfont}%
\fi\endgroup%
\begin{picture}(3600,3600)(1,-2761)
\put(251,-961){\makebox(0,0)[lb]{\smash{{\SetFigFont{12}{14.4}{\familydefault}{\mddefault}{\updefault}{\color[rgb]{0,0,0}$\left(0,\frac{1}{\log(1/p)}\right)$}%
}}}}
\put(1101,264){\makebox(0,0)[lb]{\smash{{\SetFigFont{12}{14.4}{\familydefault}{\mddefault}{\updefault}{\color[rgb]{0,0,0}$\left(\frac{\kappa}{(2dp+1)},\kappa\right)$}%
}}}}
\put(1351,-511){\makebox(0,0)[lb]{\smash{{\SetFigFont{12}{14.4}{\familydefault}{\mddefault}{\updefault}{\color[rgb]{0,0,0}$\left(\lambda,\frac{\lambda}{1-p}\right)$}%
}}}}
\put(901,-1986){\makebox(0,0)[lb]{\smash{{\SetFigFont{12}{14.4}{\familydefault}{\mddefault}{\updefault}{\color[rgb]{0,0,0}$\left(\frac{1}{\log(2d)}, \frac{1}{\log(2d)}\right)$}%
}}}}
\end{picture}%
  \caption{The set $\mathcal{D}$ in the case of the unit ball, $d=3$.}
  \label{fig:unitb}
\end{figure}

\bigskip\noindent \textbf{Proof.} (i) First consider $x$ as a function
of $y$ satisfying (\ref{basic2}). We seek the maximum, where the
derivative $x^{\prime}(y)=0$. Differentiating (\ref{basic2}) and
putting $x^{\prime}=0$, a simple calculation leads to
\[
-\log y+\log(y-x)+\log(1/p)=0,
\]
i.e.
\[
y=x/(1-p).
\]
It can be seen that this is the value of $y$ when $x$ takes its
maximum. Substituting this into (\ref{basic2}), we get
\[
x_{\max}={\frac{1}{\log(2d(1-p))}}=\lambda,
\]
verifying (\ref{xmax}).

Next consider $y$ as a function of $x$ and maximize $y$
subject to (\ref{basic2}). Again, differentiating (\ref{basic2}) with
respect to $x$ and putting $y^{\prime}=0$, we
get
\[
- \log(y-x)+\log x-\log(1/p)+\log(2d)=0
\]
from which $x=y/(1+2pd)$. Substituting in (\ref{basic2}) we get
$y_{\max}=\kappa$.

This completes the proof of Lemma 3.2(i) and the first two statements
in Lemma 3.2(ii). An easy calculation shows that if $x=0$, then
$y=1/\log(1/p)$.

Now we turn to the proof of Lemma 3.2(iii). For given $0\leq x\leq
\lambda$ consider $f(x,y)$ as a function of $y$. We have
\[
{\frac{\partial f}{\partial y}}=\log{\frac{y-x}{py}}
\]
and this is equal to zero if $y=x/(1-p)$. It is easy to see that
$f$ takes a minimum here and is decreasing if $y<x/(1-p)$ and
increasing if $y>x/(1-p)$. Moreover,
\[
{\frac{\partial^2 f}{\partial y^2}}={\frac{1}{y-x}}-{\frac{1}{y}}>0,
\]
hence $f$ is convex from below. We have for $0<x<\lambda,$ that
this minimum is
\[
f\left(x,\frac{x}{1-p}\right)=x\log((1-p)2d))=\frac{x}{\lambda}<1,
\]
and
\[
f(x,0)=x\log(2d)\left\{
\begin{array}{ll}
& < 1\,\, \mathrm{if}\,\, x<1/\log(2d), \\
& = 1\,\, \mathrm{if}\,\, x=1/\log(2d), \\
& > 1\,\, \mathrm{if}\,\, x>1/\log(2d).
\end{array}
\right.
\]
This shows that equation (\ref{basic2}) has one solution if
$0\leq x<1/\log(2d)$ and two solutions if $1/\log(2d)\leq
x<\lambda$.

For $x=\lambda$, it can be seen that $y=\lambda /(1-p)$ is the
only solution of $f(x,y)=1$.

The proof of Lemma 3.2 is complete. $\Box$

For further reference once again introduce the following notations
to describe the boundary of ${\cal D}$: for $1/\log(2d)\leq
x<\lambda$ let $y_{1,{\cal D}}(x)<y_{2,{\cal D}}(x)$ denote the
two solutions and for $0\leq x<1/\log(2d)$ let $y_{2,{\cal D}}(x)$
denote the only solution of (\ref{basic2}). Define $y_{1,{\cal
D}}(x)=x$ for $0\leq x<1/\log(2d)$ and $y_{1,{\cal
D}}(\lambda)=y_{2,{\cal D}}(\lambda)=\lambda/(1-p)$. Then the set
${\cal D}$ can be given as
$$
{\cal D}=\{0\leq x\leq \lambda,\, y_{1,{\cal D}}(x)\leq y\leq
y_{2,{\cal D}}(x)\}.
$$

For further discussion of properties of the set ${\cal D}$ see Section
6.

\renewcommand{\thesection}{\arabic{section}.}

\section{Proof of Theorems, Parts (i)}

\renewcommand{\thesection}{\arabic{section}} \setcounter{equation}{0} %
\setcounter{theorem}{0} \setcounter{lemma}{0}

In this section we prove parts (i) of the theorems in the following
order: Theorem 1.2(i), Theorem 1.1(i), Theorem 1.4(i), Theorem 1.3(i).
In the proofs the constant $c$ may vary from line to line.

\bigskip\noindent
{\bf Proof of Theorem 1.2(i).}

We say that $\mathbf{S}_j$ $(j=0,1,\ldots)$ is new (cf. \cite{DE50}) if
either $j=0$, or $j \geq 1$ and
\[
\mathbf{S}_m\neq \mathbf{S}_j,\, \, {\rm for}\, \, m=1,2,\ldots,j-1.
\]
Let $A_{j}$ be the event that $\bs_j$ is new.

Consider the reverse random walk starting from $\bs_j$, i.e.
$\bs_r':=\bs_{j-r}-\bs_j$, $r=0,1,\ldots,j$ and also the forward random
walk $\bs_r'':=\bs_{j+r}-\bs_j$, $r=0,1,2,\ldots$ Then
$\{\bs_0',\bs_1',\ldots,\bs_j'\}$ and $\{\bs_0'',\bs_1'',\ldots\}$ are
independent and so are their respective local times $\xi'$ and $\xi"$.
One can easily see that
$$
\xi(\bs_j,j)=\xi'(\boo,j)+1,\quad \xi(\bs_j+\bee_i,j)=\xi'(\bee_i,j)
$$
and
$$
\xi(\bs_j,\infty)-\xi(\bs_j,j)=\xi"(\boo,\infty),\quad
\xi(\bs_j+\bee_i,\infty)-\xi(\bs_j+\bee_i,j)=\xi"(\boo+\bee_i,\infty),
$$
hence by Lemmas 2.2 and 2.5
$$
\bp(\xi(\bs_j,\infty)=k,\xi(\bs_j+\bee_i,\infty)=\ell,\, A_j)
$$
$$
=\bp(\xi'(\boo,j)=0,\xi"(\boo,\infty)=k-1,
\xi'(\bee_i,j)+\xi"(\bee_i,\infty)=\ell)
$$
$$
=\sum_{\ell_1=0}^\ell \bp(\xi'(\boo,j)=0,\xi'(\bee_i,j)=\ell_1)
\bp(\xi"(\boo,\infty)=k-1,\xi"(\bee_i,\infty)=\ell-\ell_1)
$$
$$
\leq\sum_{\ell_1=0}^\ell\alpha^{\ell_1}
{k-1+\ell-\ell_1\choose \ell-\ell_1}\alpha^{k-1+\ell-\ell_1}
=\alpha^{k+\ell-1}\sum_{\ell_1=0}^\ell
{k-1+\ell-\ell_1\choose \ell-\ell_1}={k+\ell\choose
\ell}\alpha^{k+\ell-1}.
$$

Let $(k,\ell)\not\in((1+\varepsilon)\log n)\mathcal{B}$. Since
$g(cx,cy)=cg(x,y)$ for any $c>0$, we conclude from (\ref{joint1}) that
\[
\mathbf{P}(\xi(\mathbf{S}_j,\infty)=k,\xi(\mathbf{S}%
_j+\bee_i,\infty)=\ell, A_{j}) \leq\frac{c}{n^{1+\varepsilon}}
\]
and using this and (\ref{geo1})
\[
\mathbf{P}(\xi(\mathbf{S}_j,\infty),\, \xi(\mathbf{S}%
_j+\bee_i,\infty))\not\in ((1+\varepsilon)\log n)\mathcal{B},A_{j})
\]
\[
\leq \sum_{{(k,\ell)\not\in ((1+\varepsilon)\log n){\cal B}
\atop k\leq (1+\varepsilon)\lambda\log n}
\atop \ell\leq (1+\varepsilon)\lambda\log n}
\mathbf{P}(\xi(\mathbf{S}_j,\infty)=k,\xi(\mathbf{S}%
_j +\bee_i,\infty)=\ell,A_{j})
\]
\[
+\sum_{k>(1+\varepsilon)\lambda\log n}
\mathbf{P}(\xi(\mathbf{S}_j,\infty)=k,A_{j})
+\sum_{\ell>(1+\varepsilon)\lambda\log n}
\mathbf{P}(\xi(\mathbf{S}_j+\bee_i,\infty)=\ell,A_{j})
\]
\begin{equation}
\leq\frac{c\log^2 n}{n^{1+\varepsilon}}+
2\sum_{k>(1+\varepsilon)\lambda\log n}(1-\gamma)^k
\leq\frac{c}{n^{1+\varepsilon/2}}.
\label{epshalf}
\end{equation}

Hence selecting a subsequence $n_r=r^{4/\varepsilon}$ we have
\[
\mathbf{P}(\cup_{j\leq n_{r+1}}\cup_{i=1}^{2d} \{(\xi(\mathbf{S}_j%
,\infty), \xi(\mathbf{S}_j+\mathbf{e}%
_i,\infty))\not\in ((1+\varepsilon)\log n_r)\mathcal{B}\})
\]
\[
=\mathbf{P}(\cup_{j\leq n_{r+1}}\cup_{i=1}^{2d} \{(\xi(\mathbf{S}_j%
,\infty), \xi(\mathbf{S}_j+\mathbf{e}%
_i,\infty))\not\in ((1+\varepsilon)\log n_r)\mathcal{B}\}\cap
A_{j}) \leq
{\frac{c}{%
n_r^{\varepsilon/2}}}.
\]
 Borel-Cantelli lemma implies that with probability 1 for all
large $r$ and for all $j\leq n_{r+1}$, $i\leq 2d$ we have
$$
(\xi(\mathbf{S}_j,\infty),\xi(\mathbf{S}_j+\mathbf{e}%
_i,\infty))\in ((1+\varepsilon)\log n_r)\mathcal{B}.
$$
 It follows that with probability 1
there exists $n_0$ such that if $n\geq n_0$ then
$$
(\xi(\mathbf{S}_j,\infty),\xi(\mathbf{S}_j+\mathbf{e}%
_i,\infty))\in ((1+\varepsilon)\log n)\mathcal{B}
$$
for all $i=1,2,\dots, 2d,\, j\leq n$.

This proves (i) of Theorem 1.2. $\Box$

\bigskip
\noindent{\bf Proof of Theorem 1.1(i).}

Introduce the following notation:
\begg
\xi(\kz,(n,\infty)):=\xi(\kz,\infty)-\xi(\kz,n)
\endd

Fix $i\in\{1,2,\ldots,2d\}$ and define the following events for $j\leq
n$.

\begin{eqnarray}
B(j,n)&:=&
\{(\xi(\mathbf{S}_j,n),\xi(\mathbf{S}_j+\mathbf{e}%
_i,n))\notin ((1+\varepsilon)\log n)\mathcal{B}\},\\
B^*(j,n)&:=& \{(\xi(\mathbf{S}_j,j),\xi(\mathbf{S}_j+
\mathbf{e}_i,j))\notin ((1+\varepsilon)\log n)\mathcal{B}\},\\
C(j,n)&:=&\{\bs_m\neq \bs_j, \bs_m\neq \bs_j+\bee_i,\,
m=j+1,\ldots,n\},\\
D(j,n)&:=&
\{\xi(\bs_j,(n,\infty))>0\}\cup\{\xi(\bs_j+\bee_i,(n,\infty))>0\}.
\end{eqnarray}

Considering again the reverse random walk starting from $\bs_j$, i.e.
$\bs_r'=\bs_{j-r}-\bs_{j}$, $r=0,1,\ldots, j$ we have
$$
\xi(\bs_j,j)=\xi'(\boo,j)+1,\qquad \xi(\bs_j+\bee_i,j)=\xi'(\bee_i,j),
$$
where $\xi'$ is the local time of the random walk $\bs'$.

By (\ref{aln}) of Lemma 2.5 and (\ref{joint1}), if $(k,\ell)\notin
((1+\varepsilon)\log n)\mathcal{B}$, then
$$
\bp(\xi'(\boo,j)=k-1,\xi'(\bee_i,j)=\ell)\leq\frac{k}{k+\ell}
{k+\ell\choose k}\alpha^{k+\ell-1}\leq \frac{c}{n^{1+\varepsilon}}.
$$
Hence, as in (\ref{epshalf}), we have
$$
\bp(B^*(j,n))\leq \frac{c}{n^{1+\varepsilon/2}}.
$$

Observe that

$$B(j,n)C(j,n)D(j,n)=B^*(j,n)C(j,n)D(j.n).$$

Furthermore $\{\bs_r',\, r=0,1,\ldots,j\}$ and $\{\bs_m-\bs_j,\,
m=j,j+1,\ldots\}$ are independent. Hence

$$\bp(B(j,n)C(j,n)D(j,n))=\bp(B^*(j,n))P(C(j,n)D((j,n)).$$

Combining these with Theorem A implies
$$
\bp(B(j,n)C(j,n)D(j,n))\leq \frac
{c}{n^{1+\varepsilon/2}(n-j+1)^{d/2-1}},
$$
consequently for $d\geq 4$
\begin{equation}
\sum_{n=1}^\infty\sum_{j=1}^n
\bp(B(j,n)C(j,n)D(j,n))<\infty.
\label{conv}
\end{equation}
Hence with probability 1, there exists $n_0$ such that for $n\geq n_0$
the event $\overline B(j,n)\cup \overline C(j,n)\cup\overline D(j,n)$
occurs. We may assume that $n_0$ satisfies also the requirement in
Theorem 1.2(i). If $\overline B(j,n)$ occurs, then
$$(\xi(\mathbf{S}_j,n),\xi(\mathbf{S}_j+\mathbf{e}%
_i,n))\in ((1+\varepsilon)\log n)\mathcal{B}.$$
If $\overline D(j,n)$ occurs, then
$$
(\xi(\mathbf{S}_j,n),\xi(\mathbf{S}_j+\mathbf{e}%
_i,n))=(\xi(\mathbf{S}_j,\infty),\xi(\mathbf{S}_j+\mathbf{e}%
_i,\infty))\in ((1+\varepsilon)\log n)\mathcal{B}
$$
by Theorem 1.2(i).

Now consider $\kz\in\zd$ such that
$\xi(\kz,n)+\xi(\kz+\bee_i,n)>0$, but arbitrary otherwise and let
$L$ be the time of the last visit to $\{\kz,\kz+\bee_i\}$ before
$n$, i.e. $L:=\max\{m\leq n:\, \bs_m\in \{\kz,\kz+\bee_i\}\}$.
Then $\overline B(L,n)\cup \overline C(L,n)\cup\overline D(L,n)$
occurs for $n\geq n_0$. Since $\overline C(L,n)$ cannot occur, we
have that $\overline B(L,n)\cup\overline D(L,n)$ occurs. If
$\bs_L=\kz$, this implies

$$(\xi(\mathbf{S}_L,n),\xi(\mathbf{S}_L+\mathbf{e}%
_i,n))=(\xi(\kz,n),\xi(\kz+\bee_i,n))\in ((1+\varepsilon)\log
n)\mathcal{B}.$$
If $\bs _L =\kz+\bee_i$, then applying  the above procedure
using the unit vector $-\bee_i$ we get that
$$(\xi(\bs_L,n),\xi(\bs_L-\bee_i,n))=(\xi(\kz+\bee_i,n),\xi(\kz,n))\in
((1+\varepsilon)\log n)\mathcal{B}.$$
By symmetry of the set
$\mathcal{B}$ this implies also
$$(\xi(\kz,n),\xi(\kz+\bee_i))\in ((1+\varepsilon)\log
n)\mathcal{B}.$$

Since $i\in\{1,2,\ldots,2d\}$ is arbitrary, this
completes the proof of Theorem 1.1(i). $\Box$

\bigskip\noindent{\bf Proof of Theorem 1.4(i).}

The proof is similar to that of Theorem 1.2(i).
Let $\kz\in \zd$ and consider the unit ball centered at $\kz$. Let
now $A_j$ be the event that the random walk hits this unit ball first at
time $j$. Under this condition $(\xi(\kz,\infty),\Xi(\kz,\infty))$
has the (unconditional) distribution of
$(\xi(\boo,\infty),\Xi(\boo,\infty))$. Hence if $(k,\ell)\notin
((1+\varepsilon)\log n){\cal D}$, then by using (\ref{joint2})
$$
\bp(\xi(\kz,\infty)=k,\Xi(\kz,\infty)=\ell,A_j)\leq
\frac{c}{n^{1+\varepsilon}}.
$$
The same way as in the proof of Theorem 1.2(i) we can show the following
estimation, with the modification that whenever we have a summation by
$\ell$, $\lambda$ should be replaced by $\kappa$ and instead of using
(\ref{geo1}) we apply (\ref{ageo2}).
$$
\bp((\xi(\kz,\infty),\Xi(\kz,\infty))\notin
((1+\varepsilon)\log n){\cal D}),A_j)\leq
\frac{c}{n^{1+\varepsilon/2}}.
$$

For $n_r$ as in the proof of Theorem 1.2(i), one gets similarly
\[
\mathbf{P}(\cup_{j\leq n_{r+1}}\cup_{i=1}^{2d} \{(\xi(\mathbf{S}_j%
+\bee_i,\infty), \Xi(\mathbf{S}_j+\mathbf{e}%
_i,\infty))\not\in ((1+\varepsilon)\log n_r)\mathcal{B}\})
\]
\[
=\mathbf{P}(\cup_{j\leq n_{r+1}}\cup_{i=1}^{2d} \{(\xi(\kz%
+\bee_i,\infty), \Xi(\kz+\mathbf{e}%
_i,\infty))\not\in ((1+\varepsilon)\log n_r)\mathcal{B}\}\cap
A_{j}) \leq
{\frac{c}{%
n_r^{\varepsilon/2}}}
\]
and we can complete the proof by using Borel-Cantelli lemma.
$\Box$

\bigskip\noindent
{\bf Proof of Theorem 1.3(i).}

The proof is similar to that of Theorem 1.1(i).

Introduce the following notation:
\begg
\Xi(\kz,(n,\infty)):=\Xi(\kz,\infty)-\Xi(\kz,n).
\endd
Define $\Gamma=\Gamma_i:=\{\mathbf{e}_i+\mathcal{S}(1)\}$.
For $i\in\{1,\ldots,2d\}$ introduce, as
before, the following events for $j\leq n$.

\begin{eqnarray}
B(j,n)&:=&
\{(\xi(\mathbf{S}_j+\bee_i,n),\Xi(\mathbf{S}_j+\mathbf{e}%
_i,n))\notin ((1+\varepsilon)\log n)\mathcal{D}\},\\
B^*(j,n)&:=&
\{(\xi(\mathbf{S}_j+\bee_i,j),\Xi(\mathbf{S}_j+\mathbf{e}%
_i,j))\notin ((1+\varepsilon)\log n)\mathcal{D}\},\\
C(j,n)&:=&\{\bs_m\notin\bs_j+\Gamma,\, m=j+1,\ldots,n
\, \},\\
D(j,n)&:=&
\{\Xi(\bs_j+\bee_i,(n,\infty))>0\}.
\end{eqnarray}

Considering again the reverse random walk starting from $\bs_j$, i.e.
$\bs_r'=\bs_{j-r}-\bs_{j}$, $r=0,1,\ldots, j$ we remark
$$
\xi(\bs_j+\bee_i,j)=\xi'(\bee_i,j),\qquad
\Xi(\bs_j+\bee_i,j)=\Xi'(\bee_i,j)-1,
$$
where $\Xi'$ is the occupation time of the unit ball of the random walk
$\bs'$.

>From this we can follow the proof of Theorem 1.1(i), using
(\ref{pnn}) and (\ref{joint2}) instead of (\ref{aln}) and (\ref{joint1})
and applying Theorem 1.4(i) instead of Theorem 1.2(i). $\Box$

\renewcommand{\thesection}{\arabic{section}.}

\section{Proof of Theorems, Parts (ii)}

\renewcommand{\thesection}{\arabic{section}} \setcounter{equation}{0} %
\setcounter{theorem}{0} \setcounter{lemma}{0}

In this Section we prove parts (ii) of the Theorems.

\noindent
{\bf Proof of Theorem 1.1(ii) and Theorem 1.2(ii).}

Without loss of generality we give the proof for $i=1$. Define the
two-point set $\Upsilon:=\{\boo,\bee_1\}$. We say that $\mathbf{S}_j$
$(j=1,2,3\ldots)$ is $\Upsilon$%
-new if either $j=1$, or $j\geq 2$ and
\[
\mathbf{S}_m\notin \mathbf{S}_j+\Upsilon,\quad (m=1,2,\ldots,j-1).
\]

\begin{lemma}
Let $\zeta_n$ denote the number of $\Upsilon$-new points up to time $n$.
Then
\[
\lim_{n\to\infty}\frac{\zeta_n}{n}=1-2\alpha\quad\mathrm{a.s.}
\]
\end{lemma}

\noindent \textbf{Proof.} Define
\[
Z_j=\left\{
\begin{array}{ll}
1 & \, \mathrm{if}\,\, \mathbf{S}_j \,\, \mathrm{is}\,\, \Upsilon
{\rm-new}
\\
0 & \, \mathrm{otherwise}
\end{array}
\right.
\]
Then $\zeta_n=\sum_{j=1}^n Z_j$ and hence
\[
\mathbf{E}(\zeta_n)=\sum_{j=1}^n P(Z_j=1),
\]
\[
\mathbf{E}(\zeta_n^2)=\mathbf{E}\left(\sum_{j=1}^n\sum_{i=1}^n
Z_jZ_i\right)=
\mathbf{E}\left(\sum_{j=1}^n Z_j\right)+ 2\mathbf{E}\left(\sum_{j=1}^n%
\sum_{i=1}^{j-1}Z_jZ_i\right)
\]
\[
\leq n +2\sum_{j=1}^n\sum_{i=1}^{j-1}\bp (Z_i=1)P(Z_{j-i}=1).
\]

Considering the reverse random walk from $\mathbf{S}_i$ to
$\mathbf{S}_0=0$, we see that the event $\{Z_i=1\}$ is equivalent
to the event that the reversed random walk starting from any point
of $\Upsilon$ does not return to $\Upsilon$ up to time $i$. Using
Lemma 2.1 and 2.4 we get
\[
\bp (Z_i=1)=1-q_{\bee_1}(i)-s_{\bee_1}(i)=1-2\alpha+O(i^{1-d/2}).
\]
Hence
\[
\mathbf{E}(\zeta_n^2)\leq n+2\sum_{j=1}^n\sum_{i=1}^{j-1}
\left(1-2\alpha+O(i^{1-d/2})\right)
\left(1-2\alpha+O((j-i)^{1-d/2})\right)
\]
\[
=n(n-1)\left(1-2\alpha\right)^2+O(n^{3/2}),
\]
thus
\[
Var(\zeta_n)=O(n^{3/2}).
\]
By Chebyshev's inequality
\[
P\left(\left|\zeta_n-n\left(1-2\alpha\right)\right|>\varepsilon
n\right)\leq O\left(\frac1{\sqrt{n}}\right).
\]
Considering the subsequence $n_k=k^3$, and using Borel-Cantelli lemma and
the monotonicity of $\zeta_n$, we obtain the lemma. $\Box$

\begin{lemma}
For each $\delta>0$, there exist a subsequence $n_r$ and $r_0$ such that
if $%
r\geq r_0$ then for any $(k,\ell)\in ((1-\delta)\log n_r)\mathcal{B}
\cap%
\mathcal{Z}_d$ there exists a random integer $j_r=j_r(k,\ell)\leq
n_r$ for which
\[
(\xi(\mathbf{S}_{j_r},n_r),\xi(\mathbf{S}_{j_r}+\mathbf{e}_1,n_r))=
(\xi(\mathbf{S}_{j_r},\infty),\xi(\mathbf{S}_{j_r}+\mathbf{e}%
_1,\infty))= (k+1,\ell).
\]
\end{lemma}

\noindent\textbf{Proof.} Let $\{a_n\}$ and $\{b_n\}$ ($a_n\log n\ll b_n\ll n$%
) be two sequences to be chosen later. Define
\[
\theta_1=\min\{j>b_n:\bs_j\, {\rm is}\, \Upsilon {\rm-new}\},
\]
\[
\theta_m=\min\{j>\theta_{m-1}+b_n: \bs_j\, {\rm is} \,
\Upsilon {\rm -new},%
\}\quad m=2,3,\dots
\]
and let $\zeta_n^{\prime}$ be the number of $\theta_m$ points up
to time $n-b_n$. Obviously $\zeta_n^{\prime}(b_n+1)\geq \zeta_n$,
hence $\zeta_n^{\prime}\geq \zeta_n/(b_n+1)$ and it follows from
Lemma 5.1 that for $c<1-2\alpha$, we have with probability 1 that
$\zeta_n^{\prime}>u_n:=[cn/(b_n+1)]$ except for finitely many $n$.

For $1\leq i\leq u_n$ let
$$
\rho_0^i=0, \qquad \rho_h^i=\min\{j>\rho_{h-1}^i:
\bs_{\theta_i+j}\in \Upsilon\},\quad h=1,2,\ldots
$$

For a fixed pair of integers  $(k, \ell)$ define the following
events:
\begin{eqnarray*}
A_i&=&\{\xi(\mathbf{S}_{\theta_i},\theta_i+\rho_{k+\ell}^i)=k+1,
\xi(\mathbf{S}_{\theta_i}+\bee_1, \theta_i+\rho_{k+\ell}^i)=\ell,\\
&&\rho_h^i-\rho_{h-1}^i\leq a_n,\, h=1,\ldots,k+\ell,\,
\bs_j\not\in \bs_{\theta_i}+\Upsilon,\,
j=\theta_i+\rho_{k+\ell}^i+1,\ldots,\theta_i+b_n\},\\
B_i&=& \{\bs_j\not\in \bs_{\theta_i}+\Upsilon,\,
j>\theta_i+b_n\},\\
C_n&=&A_1B_1+\overline{A_1}A_2B_2+\overline{A_1}\overline{A_2}A_3B_3+\dots
+\overline{A_1}\dots\overline{A_{u_n-1}}A_{u_n}B_{u_n}.
\end{eqnarray*}

Note that if $(k,\ell)\in ((1-\delta)\log n){\cal B}$, then $k+\ell\leq
c\log n$ for some constant $c$, hence
$\rho_h^i-\rho_{h-1}^i\leq a_n,\, h=1,\ldots,k+\ell$ implies
$\rho_{k+\ell}^i\leq (k+\ell)a_n\leq ca_n\log n\leq b_n$ and so the
events $A_i$ are well defined and are independent, since $A_i$ depends
only on the part of random walk between $\theta_i$ and $\theta_{i+1}$.
More precisely, the events $A_1,\ldots,A_{i-1},A_iB_i$ are independent.
Moreover, $\mathbf{P}(A_i)=\mathbf{P}(A_1)$ and
$\mathbf{P}(A_iB_i)=\mathbf{P%
}(A_1B_1)$, $i=2,3,\dots$ Hence we have
\[
\mathbf{P}(C_n)=\mathbf{P}(A_1B_1)\sum_{j=0}^{u_n-1}(1-\mathbf{P}(A_1))^j=
\frac{\mathbf{P}(A_1B_1)}{\mathbf{P}(A_1)}(1-(1-\mathbf{P}(A_1))^{u_n}.
\]

\[
\mathbf{P}(\overline{C_n})\leq 1-\frac{\mathbf{P}(A_1B_1)}{\mathbf{P}(A_1)}%
+e^{-u_n\mathbf{P}(A_1)}
\]
\[
\mathbf{P}(A_1B_1)=
\mathbf{P}(D\cap\{\bs_j\not\in \Upsilon,\,
j=\rho_{k+\ell}+1,\rho_{k+\ell}+2,\ldots\})
=(1-2\alpha)\mathbf{P}(D),
\]
\begin{equation}
\mathbf{P}(A_1)= \mathbf{P}(D\cap\{\bs_j\not\in \Upsilon,\,
j=\rho_{k+\ell}+1,\ldots,b_n\})
=(1-2\alpha+O(b_n^{1-d/2}))\mathbf{P}(D), \label{pa1}
\end{equation}
where
$$
\rho_0=0, \qquad \rho_h=\min\{j>\rho_{h-1}:
\bs_{j}\in \Upsilon\},\quad h=1,2,\ldots,
$$
$$
D=\{\xi(\boo,\rho_{k+\ell})=k, \xi(\bee_1,\rho_{k+\ell})=\ell,
\rho_h-\rho_{h-1}\leq a_n,\, h=1,\ldots,k+\ell\}.
$$

In (\ref{pa1}) we used that by Lemmas 2.1, 2.4 and remembering that
$q_{\bee_1}=s_{\bee_1}=\alpha$, we have
\[
\mathbf{P}(D\cap\{\bs_j\not\in \Upsilon,\,
j=\rho_{k+\ell}+1,\ldots,b_n\})
\]
\[
\leq \mathbf{P}(D)(1-q_{\bee_1}(b_n-(k+\ell)a_n)-
s_{\bee_1}(b_n-(k+\ell)a_n))=\mathbf{P}(D)(1-2\alpha+O(b_n^{1-d/2})).
\]
Consequently,
\[
\frac{\mathbf{P}(A_1B_1)}{\mathbf{P}(A_1)}=
1+O(b_n^{1-d/2}),
\]
therefore
\[
\mathbf{P}(\overline{C_n})\leq O(b_n^{1-d/2})+e^{-cn\mathbf{P}(A_1)/b_n}.
\]

Choosing $b_n=n^{\delta/2}$, $a_n=n^{\delta/4}$, we prove
\begg
\bp(A_1)\geq\frac{c}{n^{1-\delta}}.
\label{a1}
\endd
Using (\ref{qn}) and (\ref{sn}) of Lemma 2.4 for $\kz=\bee_1$  we
get
$$
\bp(A_1)\geq (1-2\alpha){k+\ell\choose \ell}
\left(\alpha+O(a_n^{1-d/2})\right)^{k+\ell}\geq c{k+\ell\choose\ell}
\alpha^{k+\ell},
$$
since if $(k,\ell)\in {(\log n)\cal B}$, then $k+\ell=O(\log n)$.
Now (\ref{a1}) follows from Stirling formula, similarly to
(\ref{joint1}).

Using (\ref{a1}) we can verify that $\sum_r%
\mathbf{P}(\overline{C}_{n_r})<\infty$ for $n_r=r^{\rho}$ with
$\rho\delta(d-2)>4$.

By Borel-Cantelli lemma, with probability 1, $C_{n_r}$ occurs for all
but finitely many $r$. This completes the proof of Lemma 5.2. $\Box$

On choosing $\delta=\varepsilon/2$, we can see for $n_r\leq n<n_{r+1}$
\[
((1-\varepsilon)\log n)\mathcal{B}\subset ((1-\varepsilon/2)\log n_r)%
\mathcal{B}
\]
for large enough $r$ and since $\xi(\mathbf{S}_{j_r},n)$ and $\xi(%
\mathbf{S}_{j_r}+\mathbf{e}_1,n)$ do not change for $n\geq n_r$,
we have the Theorem 1.2(ii) and the first statement of Theorem
1.1(ii). The second statement in this Theorem follows by symmetry.
$\Box$

\bigskip\noindent
{\bf Proof of Theorem 1.3(ii) and Theorem 1.4(ii).}

The proof in this subsection is almost the same as in the previous
one, so we skip some details. Without loss of generality, the
proof is given for $i=1$. Let $\Gamma=\Gamma_1$ as defined in the proof
of Theorem 1.3(i), i.e. $\Gamma$ is the unit ball centered at $\bee_1$.
$\mathbf{S}_j$ $(j=1,2,3\ldots)$ is called $\Gamma$%
-new if either $j=1$, or $j\geq 2$ and
\[
\mathbf{S}_m\notin \mathbf{S}_j+\Gamma,\quad (m=1,2,\ldots,j-1).
\]

\begin{lemma}
Let $\nu_n$ denote the number of $\Gamma$-new points up to time $n$. Then
\[
\lim_{n\to\infty}\frac{\nu_n}{n}=1-p-\frac1{2d}\quad\mathrm{a.s.}
\]
\end{lemma}

\noindent \textbf{Proof.} Define
\[
Z_j=\left\{
\begin{array}{ll}
1 & \, {\rm if}\,\, \mathbf{S}_j \,\, {\rm is}\,\, \Gamma
{\rm-new}
\\
0 & \, \mathrm{otherwise}
\end{array}
\right.
\]
Then $\nu_n=\sum_{j=1}^n Z_j$.

Considering the reverse random walk from $\mathbf{S}_i$ to
$\mathbf{S}_0=0$, we see that the event $\{Z_i=1\}$ is equivalent
to the event that the reversed  random
walk starting from any point of $S(1)$ does not return to $S(1)$ up to time $%
i$. Using Lemma 2.4 we get
\[
P(Z_i=1)=1-p(i)-\frac1{2d}=1-p-\frac1{2d}+O(i^{1-d/2}).
\]

The rest of the argument is identical with that of Lemma 5.1.

\begin{lemma}
For each $\delta>0$, there exist a subsequence $n_r$ and $r_0$ such that
if $%
r\geq r_0$ then for any $(k,\ell)\in ((1-\delta)\log n_r)\mathcal{D}
\cap%
\mathcal{Z}_d$ there exists a random integer $j_r=j_r(k,\ell)\leq
n_r$ for which
\[
(\xi(\mathbf{S}_{j_r}+\mathbf{e}_1,n_r),\Xi(\mathbf{S}_{j_r}+\mathbf{e}_1,n_r))=
(\xi(\mathbf{S}_{j_r}+\mathbf{e}_1,\infty),\Xi(\mathbf{S}_{j_r}+\mathbf{e}%
_1,\infty))= (k,\ell+1).
\]
\end{lemma}

\noindent\textbf{Proof.} Let $\{a_n\}$ and $\{b_n\}$ ($a_n\log n\ll b_n\ll n$%
) be two sequences to be chosen later. Define
\[
\theta_1=\min\{j>b_n: S_j\, {\rm is}\, \Gamma  {\rm -new}\},
\]
\[
\theta_m=\min\{j>\theta_{m-1}+b_n: S_j\, {\rm is} \, \Gamma {\rm -new}%
\},\quad m=2,3,\dots
\]
and let $\nu_n^{\prime}$ be the number of $\theta_m$ points up to time
$n-b_n$. Obviously $\nu_n^{\prime}(b_n+1)\geq \nu_n$, hence
$\nu_n^{\prime}\geq\nu_n/(b_n+1)$ and it follows from Lemma 5.3 that for
$c<1-p-\frac1{2d}$, we have with probability 1 that
$\nu_n^{\prime}>u_n:=cn/(b_n+1)$ except for finitely many $n$.

Let
$$
\rho_0^i=0, \qquad \rho_h^i=\min\{j>\rho_{h-1}^i:
\bs_{\theta_i+j}\in \Gamma\},\quad h=1,2,\ldots
$$

For a fixed pair of integers $(k,\ell)$  define the following
events:
\begin{eqnarray*}
A_i&=&\{\xi(\mathbf{S}_{\theta_i}+\bee_1,\theta_i+\rho_{\ell}^i)=k,
\Xi(\mathbf{S}_{\theta_i}+\bee_1, \theta_i+\rho_{\ell}^i)=\ell+1,\\
&&\rho_h^i-\rho_{h-1}^i\leq a_n,\, h=1,\ldots,\ell,\,
\bs_j\not\in \bs_{\theta_i}+\Gamma,\,
j=\theta_i+\rho_{\ell}^i+1,\ldots,\theta_i+b_n\},\\
B_i&=& \{\bs_j\not\in \bs_{\theta_i}+\Gamma,\,
j>\theta_i+b_n\},\\
C_n&=&A_1B_1+\overline{A_1}A_2B_2+\overline{A_1}\overline{A_2}A_3B_3+\ldots
+\overline{A_1}\ldots\overline{A_{u_n-1}}A_{u_n}B_{u_n}.
\end{eqnarray*}

Similarly to the proof of Lemma 5.2,
$\mathbf{P}(A_i)=\mathbf{P}(A_1)$ and $\mathbf{P}(A_iB_i)=\mathbf{P%
}(A_1B_1)$, $i=2,3,\dots$ and
\[
\mathbf{P}(C_n)=\mathbf{P}(A_1B_1)\sum_{j=0}^{u_n-1}(1-\mathbf{P}(A_1))^j=
\frac{\mathbf{P}(A_1B_1)}{\mathbf{P}(A_1)}(1-(1-\mathbf{P}(A_1))^{u_n},
\]
\[
\mathbf{P}(\overline{C_n})\leq 1-\frac{\mathbf{P}(A_1B_1)}{\mathbf{P}(A_1)}%
+e^{-u_n\mathbf{P}(A_1)}.
\]
\[
\mathbf{P}(A_1B_1)={{\ell}\choose{k}} \left(1-p-\frac1{2d}\right)(p(a_n))^{%
\ell-k}\left(\frac1{2d}\right)^k,
\]
\[
\mathbf{P}(A_1)\leq {{\ell}\choose{k}} \left(1-p(b_n-\ell a_n)
-\frac1{2d}\right)(p(a_n))^{\ell-k}\left(\frac1{2d}\right)^k.
\]
By Lemma 2.4
\[
\frac{\mathbf{P}(A_1B_1)}{\mathbf{P}(A_1)}=1+O(b_n^{1-d/2}),
\]
therefore
\[
\mathbf{P}(\overline{C_n})\leq O(b_n^{1-d/2})+e^{-cn\mathbf{P}(A_1)/b_n}.
\]

Choosing $b_n=n^{\delta/2}$, $a_n=n^{\delta/4}$, we can prove similarly
to (\ref{a1})
\[
\mathbf{P}(A_1)\geq \frac1{n^{1-\delta}}
\]
and verify that
$\sum_r%
\mathbf{P}(\overline{C}_{n_r})<\infty$ for $n_r=r^{\rho}$ with
$\rho\delta(d-2)>4$.

By Borel-Cantelli lemma, with probability 1, $C_{n_r}$ occurs for all
but finitely many $r$. This completes the proof of Lemma 5.4. $\Box$

On choosing $\delta=\varepsilon/2$, we can see for $n_r\leq n<n_{r+1}$
\[
((1-\varepsilon)\log n)\mathcal{D}\subset ((1-\varepsilon/2)\log n_r)%
\mathcal{D}
\]
for large enough $r$ and since
$\xi(\mathbf{S}_{j_r}+\mathbf{e}_1,n)$ and $\Xi(
\mathbf{S}_{j_r}+\mathbf{e}_1,n)$ do not change for $n\geq n_r$,
we have the statements (ii) of both Theorems 1.3 and 1.4. $\Box$

\renewcommand{\thesection}{\arabic{section}.}

\section{Further discussions}

\renewcommand{\thesection}{\arabic{section}} \setcounter{equation}{0} %
\setcounter{theorem}{0} \setcounter{lemma}{0}

Observe that the following points are on the curve $g(x,y)=1$ (see
Figure 1):
$$
\left(0,\frac{1}{\log(1/\alpha)}\right),\quad
\left(\frac{1}{\log(1/\alpha)},0\right),
$$
$$
\left(\lambda,\lambda(1-\gamma)\right),\quad
\left(\lambda(1-\gamma),\lambda\right),
$$
$$
\left(\frac{1}{2\log(1/(2\alpha))},\frac{1}{2\log(1/(2\alpha))}\right).
$$

In the following discussion we are having almost sure statements,
which we will not be emphasize over and over again.

Our Theorem 1.1 shows that there are points $\kz_n$ with
$$
\xi(\kz_n,n)=0,\qquad{\rm and}\qquad \xi(\kz_n+\bee_1,n)\sim
\frac{\log
n}{\log(1/\alpha)}.
$$
On the other hand, if for a point $\kz_n$,
$$
\xi(\kz_n,n)>(1+\varepsilon)\frac{\log n}{\log(1/\alpha)},
$$
then for all of its neighbors we have $\xi(\kz_n+\bee_i,n)>c\log
n$ for some $c>0$. Moreover, if $\xi(\kz_n,n)\sim \lambda\log n$
then for all of its neighbors $\xi(\kz_n+\bee_i,n)\sim
\lambda(1-\gamma)\log n$. Roughly speaking if a point has nearly
maximal local time, it essentially determines the local time of
its neighbors, and hence the occupation time of the surface of the
unit ball around it. 

For the maximal occupation time of neighboring pairs we can obtain
$$
\lim_{n\to\infty}\frac{\sup_{\kz\in\zd}(\xi(\kz,n)+\xi(\kz+\bee_i,n))}
{\log n}=\frac1{\log\frac1{2\alpha}},
$$
and for $\kz_n$, where the sup is attained, we have, as $n\to\infty$,
$$
\xi(\kz_n,n)\sim\xi(\kz_n+\bee_i,n)\sim
\frac{\log n}{2\log\frac1{2\alpha}}.
$$

It is easy to calculate the maximal local time difference between two 
neighboring points.
$$
\lim_{n\to\infty}\frac{\sup_{\kz\in\zd}(\xi(\kz,n)-\xi(\kz+\bee_i,n))}
{\log n}=\frac1{\log\frac{1+\sqrt{1-4\alpha^2}}{2\alpha}},
$$
and for $z_n$ where the sup is attained, we have, as $n\to\infty$,

$$
\xi(\kz_n,n)\sim\frac{1+\sqrt{1-4\alpha^2}}{2\sqrt{1-4\alpha^2}}
\frac{\log n}{\log\frac{1+\sqrt{1-4\alpha^2}}{2\alpha}},
\qquad
\xi(\kz_n+\bee_i,n)\sim\frac{1-\sqrt{1-4\alpha^2}}{2\sqrt{1-4\alpha^2}}
\frac{\log n}{\log\frac{1+\sqrt{1-4\alpha^2}}{2\alpha}}.
$$

Considering now the joint behavior of the local time of a point and the
occupation time of the surface of the unit ball around it, observe that
the following points are on the curve $f(x,y)=1$ (see Figure 2):
$$
\left(0, \frac1{\log(1/p)}\right),\,
\left(\frac1{\log(2d)},\frac1{\log(2d)}\right),\,
\left(\frac{\kappa}{2dp+1},\kappa\right),\,
\left(\lambda,\frac{\lambda}{1-p}\right).
$$

As a conclusion of Theorem 1.3 we have that there are points $\kz_n$
with
$$\xi(\kz_n,n)=0 \qquad{\rm and}\qquad \Xi(\kz_n,n)\sim \frac{\log
n}{\log(1/p)}.
$$
On the other hand, if for a point $\kz_n$
$$\Xi(\kz_n,n)>(1+\varepsilon)\frac{\log n}{\log(1/p)},$$
then for its center we have $\xi(\kz_n,n)>c\log n$ for some $c>0$.
Moreover, if $\xi(\kz_n,n)\sim \lambda\log n$, then for the unit ball
$$\Xi(\kz_n,n)\sim \frac{\lambda\log n}{1-p}.$$
Roughly speaking if a point has nearly maximal local time, it
essentially determines the occupation time of the surface of the
unit ball around it.

Observe that from (\ref{pe}) it follows that
$\lambda/(1-p)=2d\lambda(1-\gamma)$, hence we may conclude that for a
ball having maximal local time at the center, the occupation time of
the surface is $2d$ times the "deterministic" local time of a
point having a neighbor with maximal local time. Consequently, all
surface points of a unit ball having maximal local time at the
center, have approximately the same local time. Moreover, if the
occupation time of the surface of a unit ball is around the
maximal value, i.e. $\Xi(\kz_n,n)\sim \kappa\log n$, then for the
local time of its center we have
$$\xi(\kz_n,n)\sim \frac{\kappa\log n}{2dp+1}.$$
Finally we conclude that even though it is natural that we can
find unit balls having the same occupation time of the surface as
the local time of its center, the fact  that it is also  possible
when this common value is fairly big is quite surprising. Namely
it is possible that
$$\xi(\kz_n,n)\sim\Xi(\kz_n,n)\sim \frac{\log n}{\log(2d)}.$$

With a  little extra computation one can easily calculate
(asymptotically) the maximal weight of the unit ball;
$$
w(\kz,n):=\xi(\kz,n)+\Xi(\kz,n),\quad w(n):=
\sup_{\kz\in \mathcal{Z}_d} (\xi(\kz,n)+\Xi(\kz,n)).
$$
This was already done in \cite{CsFRRS}. However from Theorem 1.3 we get the
following observation as well: for $d\geq 4$ if we know that
either one of the three quantities of $\xi(\kz,n),$ $\Xi(\kz,n)$
or $w(\kz,n)$ is (asymptotically) maximal, then this maximal value
uniquely determines the values of the other two (asymptotically). For
completeness here are the numerical results;

$$
\lim_{n\to \infty}\frac{w(n)}{\log n}=
-{1\over\log\left({p\over 2}+\sqrt{{p^2\over 4}+{1\over 2d}}\right)}=:
C\qquad\hspace{1cm} \mathrm{%
a.s.}  \label{last}
$$

Whenever $w(\kz_n,n)\sim C\log n$, then
$$
\xi(\kz_n,n)\sim\frac{C}{2+A}\log n,\quad  {\rm and}  \quad \Xi(\kz_n,n)\sim
C\frac{1+A}{2+A} \log n,
$$
where
$$
A=dp^2+\sqrt{d^2p^4+2dp^2}.
$$


\begin{thebibliography}{9}
\bibitem{CsFR05}  Cs\'aki, E., F\"oldes, A. and R\'ev\'esz, P.: Heavy points
of a d-dimensional simple random walk. \textit{Statist. Probab.
Lett.}, to appear.

\bibitem{CsFRRS}  Cs\'aki, E., F\"oldes, A., R\'ev\'esz, P., Rosen, J. and
Shi, Z.: Frequently visited sets for random walks. \textit{Stoch.
Process. Appl.} \textbf{115} (2005), 1503-1517.

\bibitem{DE50}  Dvoretzky, A. and Erd\H os, P.: Some problems on random walk
in space. \textit{Proc. Second Berkeley Symposium} (1951), 353--367.

\bibitem{ET60}  Erd\H os, P. and Taylor, S.J.: Some problems concerning the
structure of random walk paths. \textit{Acta Math. Acad. Sci. Hung.}
\textbf{%
11} (1960), 137--162.

\bibitem{Ham1}  Hamana, Y.: On the central limit theorem for the multiple
point range of random walk. \textit{J. Fac. Sci. Univ. Tokyo}
\textbf{39}
(1992), 339--363.

\bibitem{Ham2}  Hamana, Y.: On the multiple point range of three dimensional
random walk. \textit{Kobe J. Math.} \textbf{12} (1995), 95--122.

\bibitem{JP}  Jain, N.C. and Pruitt, W.E.: The range of transient random
walk. \textit{J. Analyse Math.} \textbf{24} (1971), 369--393.

\bibitem{Pitt}  Pitt, J.H.: Multiple points of transient random walk.
\textit{Proc. Amer. Math. Soc.} \textbf{43} (1974), 195--199.
\end{thebibliography}
\end{document}